\documentclass[11pt,twoside]{amsart}
\usepackage{amssymb}
\usepackage{dirtytalk}

\usepackage[latin1]{inputenc}
\usepackage[T1]{fontenc}
\usepackage{mathtools}
\usepackage{graphicx}
\usepackage{tikz}
\usepackage{mathabx}
\usetikzlibrary{chains}
\usepackage[OT2,T1]{fontenc}

\PassOptionsToPackage{pdfusetitle,pagebackref,colorlinks}{hyperref}
\usepackage{bookmark}
\hypersetup{
  linkcolor={red!70!black},
  citecolor={green!70!black},
  urlcolor={blue!80!black}
}

\usepackage[latin1]{inputenc}
\usepackage{amsmath}
\usepackage{amsthm}

\usepackage[all]{xy}
\usepackage{hyperref}
\newtheorem{thm}{Theorem}[section]

\newtheorem{prop}[thm]{Proposition}

\newtheorem{lem}[thm]{Lemma}
\newtheorem{cor}[thm]{Corollary}

\numberwithin{equation}{section}

\theoremstyle{definition}
\newtheorem{definition}[thm]{Definition}
\newtheorem{remark}[thm]{Remark}

\newtheorem{ex}[thm]{Example}

\DeclareFontFamily{U}{mathc}{}
\DeclareFontShape{U}{mathc}{m}{it}%
{<->s*[1.03] mathc10}{}
\DeclareMathAlphabet{\mathcal}{U}{mathc}{m}{it}
\newcommand{\kend}{\mathcal{E\mkern-3mu nd}}

\newcommand{\im}{\operatorname{im}}

\newcommand{\ch}{{\rm ch}}

\newcommand{\Aut}{{\rm Aut}}

\newcommand{\Br}{{\rm Br}}
\newcommand{\SBr}{{\rm SBr}}

\newcommand{\Coh}{{\rm Coh}}

\newcommand{\NS}{{\rm NS}}
\newcommand{\Pic}{{\rm Pic}}

\newcommand{\rk}{{\rm rk}}

\newcommand{\Spec}{{\rm Spec}}

\newcommand{\cal}{\mathcal}
\newcommand{\ka}{{\cal A}}

\newcommand{\kc}{{\cal C}}

\newcommand{\kl}{{\cal L}}

\newcommand{\km}{{\cal M}}
\newcommand{\ko}{{\cal O}}

\newcommand{\GG}{\mathbb{G}}

\newcommand{\ZZ}{\mathbb{Z}}
\newcommand{\QQ}{\mathbb{Q}}

\newcommand{\CC}{\mathbb{C}}

\newcommand{\PP}{\mathbb{P}}

\newcommand{\PPic}{\mathcal{P\mkern-3mu ic}}

\DeclareSymbolFont{cyrletters}{OT2}{wncyr}{m}{n}
\DeclareMathSymbol{\Sha}{\mathalpha}{cyrletters}{"58}

\renewcommand{\to}{\xymatrix@1@=15pt{\ar[r]&}}
\newcommand{\lto}{\xymatrix@1@=15pt{&\ar[l]}}
\renewcommand{\rightarrow}{\xymatrix@1@=15pt{\ar[r]&}}
\renewcommand{\mapsto}{\xymatrix@1@=15pt{\ar@{|->}[r]&}}
\newcommand{\mapslto}{\xymatrix@1@=15pt{&\ar@{|->}[l]&}}
\renewcommand{\twoheadrightarrow}{\xymatrix@1@=18pt{\ar@{->>}[r]&}}
\renewcommand{\hookrightarrow}{\xymatrix@1@=15pt{\ar@{^(->}[r]&}}
\newcommand{\hook}{\xymatrix@1@=15pt{\ar@{^(->}[r]&}}
\newcommand{\congpf}{\xymatrix@1@=15pt{\ar[r]^-\sim&}}
\renewcommand{\cong}{\simeq}

\newcommand{\TBC}[1]{\footnote{{\bf{\color{red} TBC $<$}} #1 \bf{\color{red} $>$ }}}
\newcommand{\Old}[1]{}
\makeatletter
\def\blfootnote{\xdef\@thefnmark{}\@footnotetext}
\makeatother

\begin{document}

\title{The Tate--{\v{S}}afarevi{\v{c}} group of a polarised K3 surface}

\author[D.\ Huybrechts, D.\ Mattei]{Daniel Huybrechts and Dominique Mattei}

\address{DH: Mathematisches Institut \& Hausdorff Center for Mathematics,
Universit{\"a}t Bonn, Endenicher Allee 60, 53115 Bonn, Germany}
\email{huybrech@math.uni-bonn.de}

\address{DM: Institute of Algebraic Geometry,
Universit{\"a}t Hannover, Welfengarten 1, 30167 Hannover, Germany}
\email{mattei@math.uni-hannover.de}

\begin{abstract} 
In \cite{HM1} we generalised the notion of the Tate--{\v{S}}afarevi{\v{c}} group $\Sha(S/\PP^1)$ of an elliptic K3 surface $S\to \PP^1$ to the Tate--{\v{S}}afarevi{\v{c}} group $\Sha(S,h)$ of a polarised K3 surface $(S,h)$.
In the present note, we complement the result by proving that for a polarised K3 surface $(S,h)$ with $h$ primitive the group $\Sha(S,h)$ parametrises bijectively all  $\Pic^0(\kc_\eta)$-torsors  that admit a `good' hyperk\"ahler compactification. Here, 
$\kc_\eta$ is  the scheme-theoretic generic curve of the linear system $|h|$.
The result is seen as the analogue of the classical fact that the Tate--{\v{S}}afarevi{\v{c}} group $\Sha(S/\PP^1)$
of an elliptic K3 surface
is the subgroup of the Weil--Ch\^atelet group of all twists that can be compactified to a K3 surface.

\end{abstract}


\maketitle
\blfootnote{The first author is supported by the ERC Synergy Grant HyperK (ID 854361).}

\section{Introduction}
Let $(S,h)$ be a polarised K3 surface over $\CC$. In particular, we assume that the polarisation $h$ is  primitive. We write $\kc\to|h|$ for the universal curve over the associated linear system and denote
its scheme geometric fibre by $\kc_\eta$, which is, therefore,  a smooth integral curve over the function field $k(\eta)$ of $|h|$. In \cite{HM1} we introduced the 
Tate--{\v{S}}afarevi{\v{c}} group $\Sha(S,h)$ of $(S,h)$ and associated with any
class $\alpha\in \Sha(S,h)$ the $\Pic^0(\kc_\eta)$-torsor $\Pic_\alpha^0(\kc_\eta)$ of $\alpha$-twisted invertible sheaves. The group $\Sha(S,h)$ was defined cohomologically by means of the special Brauer group $\SBr(S)$, 
see \cite[\S\! 4]{HM1} and \S\! \ref{sec:TSgroup} for more information.\smallskip

The classical analogue of the above situation is the case of an elliptic K3 surface $ S\to \PP^1$. Here, instead of a primitive ample $h$ one considers the pull-back of $\ko(1)$, which is primitive but of course not ample.
The Tate--{\v{S}}afarevi{\v{c}} group $\Sha(S/\PP^1)$ parametrises
by definition all K3 twists of $S\to\PP^1$ and, according to Artin--Tate \cite{Tate} and Grothendieck \cite{BrauerIII}, there exists an isomorphism $\Sha(S/\PP^1)\cong\Br(S)$.\smallskip

 In this note, we complement the cohomological definition of $\Sha(S,h)$ by the following
geometric characterisation.

\begin{thm}\label{thm:main}
Assume $(S,h)$ is a  complex K3 surface together with an ample and primitive class $h$. Then,  there exists a natural bijection
between the Tate--{\v{S}}afarevi{\v{c}} group $\Sha(S,h)$ and the 
set of isomorphism classes of $\Pic^0(\kc_\eta)$-torsors admitting a `good'
hyperk\"ahler  model
$$\xymatrix{\Sha(S,h)\ar@{<->}[r]&\{\,A=\Pic^0(\kc_\eta)\text{\rm-torsor}\mid \exists~ A\subset X \text{ \rm good HK model}\,\}.}$$
\end{thm}

By definition, a good hyperk\"ahler model compactifies $A\to\eta$ to a Lagrangian fibration
 $X\to |h|$ with all fibres required to
have at least one generically reduced component, see \S\! \ref{subsec:good} for the precise definition and further details.
 Quite possibly,
under our assumptions, every torsor admitting a compactification by some smooth projective
hyperk\"ahler variety at all also admits a hyperk\"ahler compactification that is good in our sense.

\smallskip

The result has two parts (see also \S\! \ref{sec:Proof}):
\begin{enumerate}
\item[(i)] The $\Pic^0(\kc_\eta)$-torsors parametrised by $\Sha(S,h)$ are pairwise non-isomorphic.
This will be addressed in \S\! \ref{sec:injec}.
\item[(ii)] Among all $\Pic^0(\kc_\eta)$-torsors, one obtains exactly those that admit a geometrically satisfying compactification, see (\ref{eqn2}) below for the well-known analogue of elliptic K3 surfaces.
This part will be dealt with in \S\! \ref{sec:Image}.
\end{enumerate}
\medskip

The injectivity statement (i) relies on the techniques developed in \cite{HM1} (and could have been proved there already), but for the the surjectivity
in (ii) we make extensive use of recent work of Kim \cite{Kim}.\smallskip

\noindent
{\bf Acknowledgment.} The first author wishes to thank the organisers and the participants of the
program `Arithmetic geometry of K3 surfaces' at the Bernoulli Center for the invitation and for questions
and comments.  We wish to thank the anonymous referee of \cite{HM1} whose numerous comments and constructive suggestions have been invaluable also for writing this paper. We are most grateful to Yoonjoo Kim for helpful comments on a first version of the present note.

\section{Review: Elliptic K3 surfaces}
In order to motivate our results and to increase readability, we briefly recall the notion of the Tate--{\v{S}}afarevi{\v{c}}
group of an elliptic K3 surface, emphasising the modular perspective.\smallskip

\subsection{Classical theory} We start with an elliptic K3 surface $S_0\to \PP^1$ over an algebraically closed field $k$, which one might as well assume to be  $\CC$, as we will apply Hodge theory eventually.
 Note that by definition, an elliptic K3 surface comes with a section. If we only assume that the smooth fibres are isomorphic to elliptic curves without fixing the origin by means of a section then we will speak of genus one fibred K3 surfaces.\smallskip

A \emph{twist} of $S_0\to\PP^1$ consists of a K3 surface $S$ together with a
genus one fibration $S\to\PP^1$ and an open immersion $\Pic^0(S/\PP^1_{\rm sm})\,\hookrightarrow S_0$ of the relative Jacobian of its smooth part
$S\to\PP^1_{\rm sm}$. In other words, we
fix an isomorphism of the minimal smooth compactification
of the relative Jacobian $\Pic^0(S/\PP^1_{\rm sm})$  with $S_0$ relative over $\PP^1$.\smallskip

By definition, the Tate--{\v{S}}afarevi{\v{c}} group is the group of all twists:
$$\Sha(S_0/\PP^1)\coloneqq\{\,\text{twists}\,\}/_\cong.$$
At this point, it is simply a set of isomorphism classes. Its group structure
becomes apparent by observing that any twist $S\to \PP^1$
of $S_0\to \PP^1$ defines a torsor of the generic fibre $E_0=S_{0\eta}$.
More precisely, $\Sha(S_0/\PP^1)$ can be realised as a subgroup
of the Weil--Ch\^atelet group of $E_0$:
$$\Sha(S_0/\PP^1)\,\hookrightarrow \text{WC}(E_0)\cong H^1(\eta,E_0).$$

It is known that the fibres of $S_0\to \PP^1$ and of any twist $S\to \PP^1$ over
any point $t\in \PP^1$ are isomorphic, cf.\ \cite[Lem.\ 5.2]{FM}. In fact, it is this condition that describes
$\Sha(S_0/\PP^1)$ as a subgroup of the much larger $\text{WC}(E_0)$. 

Summarising the discussion so far, there is a natural identification
\begin{equation}\label{eqn2}
\xymatrix{\Sha(S_0/\PP^1)\ar@{<->}[r]&\{\,E=E_0\text{-torsor }\mid \exists~ E\subset S \text{ K3}\, \}.}
\end{equation}

This geometric definition of $\Sha(S_0/\PP^1)$ is accompanied by a global
cohomological description in terms of the Brauer group. Recall
that the Brauer group of equivalence classes of Azumaya algebras
can be identified with the cohomological Brauer group, i.e.\ there exists an isomorphism
$$\Br(S_0)\coloneqq\{\,\ka=\text{Azumaya}\,\}/_\sim\cong H^2(S_0,\GG_m).$$

According to Artin--Tate \cite{Tate} and Grothendieck \cite{BrauerIII}, there also exists a natural isomorphism $\Br(S_0)\cong {\Sha(S_0/\PP^1)}$. The proof is purely
cohomological and uses the Hochschild--Serre spectral sequence
$$\xymatrix{\Br(S_0)\ar@{^(->}[d]\ar@{-}[r]^-\cong&\Sha(S_0/\PP^1)\ar@{^(->}[dd]\\
\Br(E_0)\ar[d]_-\cong&\\
H^1(\eta,\Pic(E_0))\ar@{-}[r]^-\cong& H^1(\eta,E_0).}$$

\subsection{Via moduli spaces}
There is a more geometric argument to prove the isomorphism $\Br(S_0)\cong \Sha(S_0/\PP^1)$ 
that relies on moduli spaces of (twisted) sheaves. We will give
both directions of this isomorphism. It is not too difficult to prove that they are inverse to each other
(possibly with a sign) and that they coincide with the isomorphism by Artin--Tate and Grothendieck.\smallskip


$\bullet$ $\pmb{\Sha(S_0/\PP^1)\to\Br(S_0)}$: Fix a twist $S\in \Sha(S_0/\PP^1)$.
Then by definition $S_0$ is isomorphic to the minimal smooth compactification
of the relative Jacobian $\Pic^0(S/\PP^1_{\rm sm})$ of $S\to \PP^1_{\rm sm}$. 
A modular compactification is provided by the moduli space $M_S(0,f,0)$ of stable sheaves with Mukai vector $(0,f,0)$, i.e.\ of stable sheaves (with respect to an appropriate polarisation) with  numerical invariants $(\rk,{\rm c}_1,\chi)=(0,f,0)$,
which is the Mukai vector of any invertible sheaf of degree zero on a smooth fibre
of $S\to\PP^1$. The situation is summarised by
$$S_0\cong\overline\Pic^0(S/\PP^1_{\rm sm})\cong M_S(0,f,0).$$

Observe that the moduli space $M_S(0,f,0)$
is typically not fine, i.e.\ there does not exist a universal sheaf on $S\times S_0\cong S\times M_S(0,f,0)$. The obstruction to the existence of a universal
sheaf is a class $\alpha\in H^2(S_0,\GG_m)\cong \Br(S_0)$. Then mapping a twist
$S\to\PP^1$ to this obstruction class defines a map $$\Sha(S_0/\PP^1)\to \Br(S_0)\text{, }S\mapsto \alpha.$$


\smallskip

Note that the moduli space perspective also allows one to show (again) that at least
all smooth fibres of $S_0\to \PP^1$ and of a twist $S\to \PP^1$ 
over $t\in \PP^1$ are isomorphic to each other. So  \'etale locally or analytically locally both families,
or at least their regular parts, are identical, the difference is in the gluing.
\smallskip

$\bullet$ $\pmb{\Br(S_0)\to \Sha(S_0/\PP^1)}$: We start with a Brauer class $\alpha\in \Br(S_0)$  and consider the moduli space of stable $\alpha$-twisted
sheaves $S_\alpha\coloneqq M_{S_0,\alpha}(0,f,\ast)$. Here, the numerical
invariants are again such that the sheaves are generically invertible sheaves on
the fibres (this corresponds to the first two entries), but we are intentionally vague about the last one `$\ast$'.\footnote{We cannot talk about $\chi(F)$ of a
twisted sheaf $F$. This seems like a technical point, but it is crucial for our approach to the Tate--{\v{S}}afarevi{\v{c}} group. We will come back to this briefly in \S\! \ref{sec:IntroSBr}, see also \cite[\S\! 2]{HM1}.} This defines a map $$\Br(S_0)\to \Sha(S_0/\PP^1)\text{, } \alpha\mapsto S_\alpha,$$
which is inverse to the above.


\section{Tate--{\v{S}}afarevi{\v{c}} group of a polarised K3 surface}\label{sec:TSgroup}
We recall the basic features of the  Tate--{\v{S}}afarevi{\v{c}} group of a polarised K3 surface $(S,h)$
and explain what it is expected to parametrise geometrically.
\subsection{Cohomological approach}\label{sec:IntroSBr}
In \cite{HM1} we introduced the Tate--{\v{S}}afarevi{\v{c}} $\Sha(S,h)$
of a polarised K3 surface $(S,h)$ and proved that there exists
a short exact sequence
\begin{equation}
0\to\ZZ/m\ZZ\to\Sha(S,h)\to \Br(S)\to0
\end{equation}
linking it to the Brauer group of $S$. Here, the positive integer $m$ is the divisibility of $h$, i.e.\ $m$ is determined by $m\ZZ=(\NS(S).h)$.
The group $\Sha(S,h)$ was introduced as a certain quotient of the subgroup
$\SBr(S,h)\subset\SBr(S)\cong H^2(S,\QQ/\ZZ)$ of all classes  $\alpha\in\SBr(S)$
for which the twisted Picard scheme $\Pic_\alpha^0(\kc_\eta)$ of the scheme-theoretic
generic fibre $\kc_\eta$ of the universal family $\kc\to|h|$ is not empty. Note that if  $\alpha$
is fixed as a class in $\SBr(S)$ and not merely in its quotient $\Br(S)$, then one can talk about the
degree of an $\alpha$-twisted invertible sheaf on $\kc_\eta$ and, therefore, $\Pic^0_\alpha(\kc_\eta)$ is well defined, cf.\ \cite[\S\! 2]{HM1}.\smallskip

In \cite[\S\! 4]{HM1} we proved that $\SBr(S,h)$ is the annihilator of $h\in \NS(S)$ for the intersection
pairing $\SBr(S)\times \NS(S)\to \QQ/\ZZ$. The description of 
$\SBr(S)$ as $H^2(S,\QQ/\ZZ)$ suggests to think of 
$\SBr(S,h)$ as the primitive cohomology\footnote{Thanks to Emma Brakkee and
Tony V\'arilly-Alvarado for this observation.}$$\SBr(S,h)\cong H^2(S,\QQ/\ZZ)_{\rm pr}.$$
In particular, a class
$\gamma=(1/r)L\in \NS(S)\otimes\QQ/\ZZ\subset \SBr(S)$ is contained in $\SBr(S,h)\subset\SBr(S)$ if and only if
$(h.\gamma)\in\ZZ$, cf.\ \cite[Prop.\ 4.2]{HM1}.\smallskip

This smaller group $\SBr(S,h)\subset\SBr(S)$ still surjects onto $\Br(S)$ so that there exists a short exact
sequence $$0\to A(S,h)\to \SBr(S,h)\to \Br(S)\to 0.$$
 Then, the group $\Sha(S,h)$ is defined as the quotient of $\SBr(S,h)$ by the kernel
of the natural map $A(S,h)\to \ZZ/m\ZZ$ given by the intersection with $h$.\smallskip

As an aside, we observe that if $h$ generates $\NS(S)$, then $A(S,h)$ is isomorphic to the discriminant $A(S)=\NS(S)^\ast/\,\NS(S)$.

\subsection{Via Weil--Ch\^atelet}
On the other hand, from a geometric point of view and in analogy to the case of elliptic K3 surfaces, we want the Tate--{\v{S}}afarevi{\v{c}} group  of the polarised K3 surface $(S,h)$ to be isomorphic to a certain
subgroup of the Weil--Ch\^atelet group 
$$``\Sha(S,h)" \subset H^1(\eta,\Pic^0(\kc_\eta))
$$ of the Picard group of the generic fibre $\kc_\eta$, which paramatrises particularly `nice' $\Pic^0(\kc_\eta)$-torsors.
We will come back to this in \S\! \ref{sec:Image}.\smallskip

Note that also the Brauer group of $S$ can be linked to the Weil--Ch\^atelet group of the full Picard
scheme $\Pic(\kc_\eta)$ (use again the Hochschild--Serre spectral sequence):
$$\xymatrix{\Br(S)\cong\Br(\kc)\ar@{^(->}[r]&\Br(\kc_\eta)\ar[r]& H^1(\eta,\Pic(\kc_\eta)).}$$
We take this as an encouraging sign.
 However, in order to restrict to $\Pic^0(\kc_\eta)$-torsors, one
needs to work with the special Brauer group. This is related to the fact that
the inclusion $\Pic^0(\kc_\eta)\,\hookrightarrow \Pic(\kc_\eta)$ induces  a surjection $H^1(\eta,\Pic^0(\kc_\eta))\twoheadrightarrow H^1(\eta,\Pic(\kc_\eta))$ which typically is not an isomorphism, for its kernel is $\ZZ/m\ZZ$.

\begin{remark} Note that by assumption, an elliptic K3 surface
$S_0\to \PP^1$ comes with a section. There is nothing this 
could be replaced by for a  general polarised K3 surface $(S,h)$. In this sense, the higher-dimensional
situation is rather analogous to the case of genus one fibred
K3 surfaces $S\to\PP^1$ without a section. Although,  classically there is no Tate--{\v{S}}afarevi{\v{c}} group
of a genus one fibred K3 surface without a section, one can use the usual short exact sequence
\begin{equation}\label{eqn1}
\xymatrix{0\ar[r]&\langle \alpha\rangle\ar[r]& \Br(S_0)\cong\Sha(S_0/\PP^1)\ar[r]& \Br(S)\ar[r]&0}
\end{equation}
and define the $\Sha$ of $S\to \PP^1$ without 
a section simply as the $\Sha$ of its Jacobian fibration $S_0\to\PP^1$, see also \cite[\S\! 5.2]{HM1}.
However,
this quick fix is not available for a polarised K3 surface $(S,h)$. There is just
no replacement for the relative Jacobian which would again be a K3 surface.
\end{remark}

\Old{
\begin{remark} Ideally, we would like the condition on $(S,h)$ to be generic just to mean that all curves $C\in |h|$ are integral. This does hold whenever $\Pic(S)=\ZZ h$ but describes a dense Zariski open subset of all polarised K3 surfaces. Unfortunately, the argument above is less precise, so that we currently cannot identify exactly the Zariski open set for which the result holds.
We only proved that it exists.
\end{remark}}

\section{Injectivity}\label{sec:injec}

In this section we prove that no two twisted Picard schemes $\Pic_{\alpha_1}^0(\kc_\eta)$  and
$\Pic_{\alpha_2}^0(\kc_\eta)$ are isomorphic as  $\Pic^0(\kc_\eta)$-torsors 
unless $\alpha_1=\alpha_2$ as classes in $\Sha(S,h)$.

\subsection{Special Brauer classes inducing trivial torsors} The main step in the proof of the injectivity
in Theorem \ref{thm:main} is the following.

\begin{prop}\label{prop:Inj}
For any polarised K3 surfaces $(S,h)$ with $h$ primitive, the
kernel of the map $$\SBr(S,h)\to H^1(\eta,\Pic^0(\kc_\eta)), ~~
\alpha\mapsto \Pic_\alpha^0(\kc_\eta)$$
coincides with the subgroup of all elements $\gamma\in\SBr(S,h)\cap (\NS(S)\otimes\QQ/\ZZ)$
with $(h.\gamma)\in m\ZZ$.
\end{prop}


\begin{proof} One direction was already proved in \cite[Prop.\ 4.5]{HM1}. For completeness,
we briefly recall the argument.
Assume $\alpha\in \SBr(S,h)$ is a class $\gamma=(1/r)\cdot L\in \NS(S)\otimes\QQ/\ZZ$,
with $L\in \NS(S)$, which satisfies $(h.\gamma)\in m\ZZ$. Then there exists an $L'\in \NS(S)$ with $(h.L')+(1/r)(h.L)=0$. Hence, replacing $L$ by $L+rL'$, we may assume that $(h.L)=0$. If now $F$ is any locally free sheaf on $S$ of rank $r$ and with
$\det(F)\cong L$, then  $E\mapsto E\otimes F|_{\kc_\eta}$ defines an isomorphism
$\Pic^0(\kc_\eta)\congpf \Pic_\alpha^0(\kc_\eta)$ of $\Pic^0(\kc_\eta)$-torsors. Thus,
$\Pic_\alpha^0(\kc_\eta)$ is trivial.\smallskip

The converse is more complicated and we proceed in two steps. The first step uses the transcendental lattice to reduce to the case $\alpha\in \NS(S)\otimes\QQ/\ZZ$ while the second relies on a more
geometric argument. An alternative proof, which uses the N\'eron--Severi lattice and assumes the K3 surface to have Picard number one, will be sketched in Remark \ref{rem:altproof}.
\smallskip

(i) Assume $\Pic_\alpha^0(\kc_\eta)$ is the trivial torsor. \emph{Claim:} Then $\alpha\in \SBr(S,h)\cap (\NS(S)\otimes\QQ/\ZZ)$, i.e.\ $\alpha$ becomes trivial under the projection
$\SBr(S,h)\twoheadrightarrow\Br(S)$, $\alpha\mapsto\bar\alpha$.
\smallskip

Since $(S,h)$ is a polarised K3 surface with $h$ primitive, there exists a smooth projective hyperk\"ahler compactification $M_\alpha$ of $\Pic_\alpha^0(\kc_\eta)$. It is  provided by the moduli space of stable $\alpha$-twisted sheaves that are numerically equivalent to 
any degree zero $\alpha$-twisted invertible sheaf on a smooth curve in $|h|$:
\begin{equation}\label{eqn:MS}
\Pic_\alpha^0(\kc_\eta)\subset\overline\Pic_\alpha^0(\kc_\eta)=M_\alpha.
\end{equation}
Here, stability is taken with respect to a generic polarisation. (The choice of the polarisation
also determines how the degrees of a stable sheaf with support $C\in |h|$ are distributed over the irreducible components
of  $C$. However, for our purposes these subtleties are irrelevant.)
In the generic case that all curves $C\in |h|$ are integral, for example when $\rho(S)=1$,  semi-stability implies stability for every polarisation.  
For the special case of trivial  $\alpha$ we introduce the notation
\begin{equation}\label{eqn:MS0}\Pic^0(\kc_\eta)\subset\overline\Pic^0(\kc_\eta)=M,
\end{equation}
where $M$ is the moduli space of stable sheaves with Mukai vector $v_0=(0,h,-h^2/2)$.\smallskip

If $\Pic_\alpha^0(\kc_\eta)$ is a trivial torsor, then $M_\alpha$ and $M$ are birational hyperk\"ahler manifolds. In this case, there exists a Hodge isometry $H^2(M_\alpha,\ZZ)\cong H^2(M,\ZZ)$ between
their weight two Hodge structures and, as a consequence, a Hodge isometry $T(M_\alpha)\cong T(M)$ between their transcendental lattices. By a result of O'Grady \cite[Main Thm.]{OG}, there also
exists a Hodge isometry  $H^2(M,\ZZ)\cong v_0^\perp\subset\widetilde H(S,\ZZ)$, which 
yields a Hodge isometry $T(M)\cong T(S)$ (despite, potentially, the absence of a universal family). The analogous result in the twisted case is due to Yoshioka \cite[Thm.\ 3.19]{Yosh} and provides us with a Hodge isometry
$T(M_\alpha)\cong T(S,\bar\alpha)$, where the latter is the kernel of $\bar\alpha\colon T(S)\to\QQ/\ZZ$,
cf.\ \cite{Cal,HS}.
Altogether, we find a (Hodge) isometry $$T(S,\bar\alpha)\cong T(M_\alpha)\cong T(M)\cong T(S).$$
By comparing discriminants, this immediately shows that $\bar\alpha$ is the trivial Brauer class.
\smallskip

(ii) Assume now that $\alpha\in\SBr(S,h)$ is contained in $\NS(S)\otimes\QQ/\ZZ$, i.e.\ it is given by a class $\gamma=(1/r)L$. \emph{Claim:} Then $(h.\gamma)\in m\ZZ$.
\smallskip

We shall give an argument that avoids all lattice theory, but uses the existence of a modular Lagrangian compactification. Under the additional assumption that $\Pic(S)$ is generated by $h$ a more direct
argument exists, see Remark \ref{rem:altproof} below. \smallskip

By assumption, we have $\Pic_\alpha^0(\kc_\eta)\cong \Pic^d(\kc_\eta)$ with $d=(h.\gamma)\in \ZZ$, see \cite[Ex.\ 4.10]{HM1}, which can be compactified by the moduli space $M(v_d)$ of stable sheaves with Mukai vector $v_d=(0,h,d-h^2/2)$:
$$\Pic_\alpha^0(\kc_\eta)\cong\Pic^d(\kc_\eta)\subset \overline\Pic^d(\kc_\eta)=M(v_d).$$

Assume now that $\Pic_\alpha^0(\kc_\eta)\cong\Pic^d(\kc_\eta)$ is trivial and fix a  trivialisation $\Pic^d(\kc_\eta)\cong\Pic^0(\kc_\eta)$. Then
the origin of $\Pic^0(\kc_\eta)$, corresponding to the trivial invertible sheaf on $\kc_\eta$, defines a $k(\eta)$-rational point of $\Pic^d(\kc_\eta)$.  Its closure $\Sigma\subset M(v_d)$ defines a rational
section of the projection $\pi\colon M(v_d)\to|h|$ and we pick a desingularisation $\tilde\Sigma\to\Sigma$.
Note that if all $C\in |h|$ and hence all fibres of $\pi$ are integral, then $\Sigma\subset M(v_d)$ is actually a section, cf.\  \cite[\S\! 6.4.1]{Kim}. \smallskip

The moduli space $M(v_d)$ will in general not be fine,  but there always exists a Brauer class
$\beta\in\Br(M(v_d))$ such that a universal family exists as a $(\beta,1)$-twisted sheaf on
$M(v_d)\times_{|h|}\kc$. As $\tilde\Sigma$ is smooth, projective, and rational,  its Brauer group
$\Br(\Sigma)$ is trivial. Hence, the pull-back of the twisted universal family to $\tilde\Sigma\times_{|h|}\kc$, which is birational to $\kc$ defines an untwisted invertible sheaf of  relative degree $d$ on the generic fibre $\kc_\eta$. \smallskip

As explained in \cite[\S\! 4]{HM1}, an invertible sheaf of degree $d$ on the generic fibre $\kc_\eta$ can be spread out to eventually define an invertible sheaf  $L$ on $S$ with $(L.h)=d$. 
In other words, since $\kc$ is smooth, the restriction map
$\Pic(\kc)\twoheadrightarrow \Pic(\kc_\eta)$ is surjective. By definition of $m$, this immediately implies $d\in m\ZZ$.
\end{proof}

\begin{remark}\label{rem:altproof}
Let us explain an alternative and more direct argument for (ii) in the proof above under the assumption
that $\Pic(S)$ is generated by $h$. As before,  we  write $\Pic_\alpha^0(\kc_\eta)\cong \Pic^d(\kc_\eta)$, and observe that  by tensoring with the appropriate invertible sheaf we may assume $0\leq d< m=h^2$.
With $M(v_d)$ as before, 
 there is a chain of Hodge isometries resp.\ Hodge isometric embeddings
$$\widetilde H(S,\ZZ)\supset v_d^\perp\cong H^2(M(v_d),\ZZ)\cong H^2(M_\alpha,\ZZ)\cong H^2(M,\ZZ)\cong v_0^\perp\subset \widetilde H(S,\ZZ),$$
which on the algebraic parts leads to an isometry
$(v_d^\perp)^{1,1}\cong (v_0^\perp)^{1,1}$.\smallskip

The proof then concludes by the following lattice argument.  Observe that $${\rm disc}(v_0^\perp)^{1,1}=4 \text{ and }{\rm disc}(v_d^\perp)^{1,1}=(h^2)^2/{\rm gcd}^2,$$
where ${\rm gcd}\coloneqq{\rm gcd}(h^2,d-h^2/2)$. For example, 
$(v_0^\perp)^{1,1}$ is integrally generated by the two vectors $(0,0,1),(2,-h,0)\in \widetilde H^{1,1}(S,\ZZ)$ with the intersection matrix $\left(\begin{matrix}0&2\\2&h^2\end{matrix}\right)$. Hence, $h^2/2={\rm gcd}$ and, therefore, $d=0$, i.e.\ $(h.\gamma)\in m\ZZ$. \smallskip
\end{remark}

\begin{remark}\label{rem:Picstable}
For later use we note that any invertible sheaf $L\in \Pic(C)$ on a possibly non-reduced or reducible
curve  $C=\sum m_iC_i\in|h|$ with $\deg(L|_{C_i})=0$  for each integral  component $C_i$
 is actually stable.  Note that at this point we do not have to assume $h$ primitive. Since the assertion is purely numerical, it suffices to prove the
 assertion for the structure sheaf $\ko_C$ itself. The usual
 stability inequality, cf.\ \cite[\S\! 1.2]{HL}, for a pure quotient of $\ko_C$, which is necessarily of the form $\ko_{C'}$  of some $C'\subset C$, reads $(C'.C')\leq (h.C')$ which is the 1-connectedness
 of big and nef divisors on K3 surfaces, cf.\ \cite[Ch.\ 2, Rem.\ 1.5]{HuyK3}.\smallskip
 
 Later we will
 phrase this as an open inclusion $\Pic^0(\kc/|h|)\,\hookrightarrow M$, see \S\! \ref{sec:UntwiNeron}.
 Observe that the argument shows stability of $\ko_C$ with respect to the fixed polarisation $h$, which might
 not be generic with respect to the Mukai vector $v(\ko_C)=(0,h,-h^2/2)$. However, stability of a fixed sheaf is an  open condition on the ample class and, therefore, $\ko_C$  remains stable for any (and hence a generic) polarisation close to $h$.\smallskip
 
 Note that line bundles of non-zero degree on components of $C$
 may not be stable, cf.\ \cite{dCRS,IGB,Kim}.
\end{remark}

\subsection{Comparison and injectivity statement}

In \cite[\S\! 4]{HM1} we explained that the more direct approach leads to two obstruction classes
in the Brauer group $\Br(k(\eta))$ of the function field of $|h|$ (which is, of course, a non-trivial group). The arguments above use the global
structure of the hyperk\"ahler model of $\Pic_\alpha^0(\kc_\eta)$ provided by the moduli space $M_\alpha$ of stable $\alpha$-twisted sheaves.\smallskip

As by definition $\Sha(S,h)$ is the quotient of $\SBr(S,h)$
by the subgroup $\{\,\gamma\mid (h.\gamma)\in m\ZZ\,\}$, the proposition  immediately leads to the
injectivity claim in Theorem \ref{thm:main}, which we state again as follows.

\begin{cor}\label{cor:Inj} For any polarised K3 surface $(S,h)$ with $h$ primitive, the map
\begin{equation}\label{eqn:inj}
\Sha(S,h)\, \hookrightarrow H^1(\eta,\Pic^0(\kc_\eta)),~\alpha\mapsto
\Pic_\alpha^0(\kc_\eta)
\end{equation} is injective.\qed
\end{cor}

\section{The image}\label{sec:Image}
We are now addressing the second question: How do we characterise those
$\Pic^0(\kc_\eta)$-torsors that are of the form $\Pic^0_\alpha(\kc_\eta)$ for
some $\alpha\in\Sha(S,h)$? Note that from the classical case of elliptic K3 surfaces, we
know already that not every $\Pic^0(\kc_\eta)$-torsor can be of the form $\Pic^0_\alpha(\kc_\eta)$. 
 In other words, (\ref{eqn:inj}) is certainly not surjective

\subsection{Good hyperk\"ahler models}\label{subsec:good}
We consider a projective space $B=\PP^n$ with its scheme-theoretic generic
point $\eta\in B$ and an abelian variety $A_0$ over the function field $k(\eta)$.
\begin{definition}\label{def:goodHK}
A \emph{good hyperk\"ahler model} of an $A_0$-torsor $A$ consists of a smooth projective hyperk\"ahler variety $X$ together with a Lagrangian fibration $X\to B$ such that:
\begin{enumerate}
\item[(i)] The generic fibre $X_\eta$ is isomorphic to $A$ and 
\item[(ii)] Every closed fibre $X_t$, $t\in B$, contains a smooth point.
\end{enumerate}
\end{definition}
Condition (ii) is of course equivalent to saying that $X_t$ has at least one generically reduced irreducible
component. According to \cite[Rem.\ 1.3]{Campana}, this condition is equivalent to $X_t$ not being multiple or, more precisely, the closed fibre $X_t$ of $X_\Delta$ over a generic disk $t\in\Delta\subset B$ cannot be written as $X_t=m\cdot F$ with $m>1$. For the scheme-theoretic notion of multiplicity
see \cite[TAG 02QS]{Stacks}.
\begin{ex}
Note that any genus one fibration $S\to \PP^1$ of a K3 surface  satisfies (ii). There might exist fibres that 
are not reduced, but no fibre is everywhere non-reduced, i.e.\ every fibre admits at least one (not necessarily unique) component that  is reduced, cf.\ \cite[Ch.\ 11]{HuyK3}.
\end{ex}

Assume now that $(S,h)$ is a polarised K3 surface which always includes $h$ primitive. As before,
we denote the universal family by $\kc\to|h|$ and its scheme-theoretic generic fibre by $\kc_\eta$. As the main example of good hyperk\"ahler models, we have the following.

\begin{lem}\label{lem:GoodHK}
For any $\alpha\in\Sha(S,h)$ the $\Pic^0(\kc_\eta)$-torsor $\Pic^0_\alpha(\kc_\eta)$ admits a good hyperk\"ahler model.
\end{lem}

\begin{proof}
A good hyperk\"ahler model is provided by the moduli space
$M_\alpha$ of stable $\alpha$-twisted sheaves, see (\ref{eqn:MS}). Indeed, the generic
fibre of $M_\alpha\to|h|$ gives back $\Pic^0_\alpha(\kc_\eta)$ by definition, so (i) holds true.\smallskip

In order to verify (ii), we have to show that the fibre $M_{\alpha t}$ of $M_\alpha\to|h|$ over a point $t\in |h|$ corresponding to   a possibly singular, reducible, non-reduced curve $C=\kc_t\subset S$ contains a smooth point. Recall from Remark \ref{rem:Picstable} that for $\alpha$ trivial, the fibre
contains $\Pic^0(C)$ which is smooth. Ideally, one would like to argue that $M_{\alpha t}$ contains the smooth $\Pic^0_{\alpha|_C}(C)$. There are two problems:
\begin{enumerate}
\item[$\bullet$] Is $\Pic^0_{\alpha|_C}(C)$ non-empty? 
\item[$\bullet$] Is $\Pic^0_{\alpha|_C}(C)$ contained in $M_{\alpha t}$?
\end{enumerate}
Since $\bar\alpha|_C$ is trivial, cf.\ \cite[Prop.\ 3.6.6]{CTS}, $\Pic^0_{\alpha|_C}(C)\cong\Pic^d(C)$ for some $d$ and the latter is non-empty if $C$ is not multiple. Indeed, if $C=\sum m_iC_i$ is not multiple, which always holds under the assumption that $h$ is primitive, then ${\rm gcd}(m_i)=1$. But then there exists
an invertible sheaf $L$ on $C$ with $\sum m_i\deg(L|_{C_i})=d$ and so $\Pic^d(C)\ne\emptyset$.
For the second point, even the existence of a single (twisted) line bundle on $C$ of the correct degree that is stable
as a sheaf on $S$ is not clear. In particular, we do not actually know whether $\Pic^0_\alpha(\kc/|h|)$ is contained
in $M_\alpha$. Note, however, that over the smooth locus we do have an open
embedding, i.e.\ $\Pic^0_\alpha(\kc/|h|_{\rm sm})\subset M_\alpha$.\smallskip

So, to address the second problem above, we let $M_\alpha^{\rm simp}$ be the moduli space of simple twisted sheaves, which is 
only a non-separated algebraic space. Then, of course, $\Pic^0_\alpha(\kc/|h|)\subset M_\alpha^{\rm simp}$. Indeed, for fixed $t\in |h|$ and $C\coloneqq\kc_t$ we pick a $L\in \Pic_{\alpha|_C}^0(C)$ and view it as a smooth
point of the fibre of $M_\alpha^{\rm simp}\to|h|$ over $t$. Next, choose a generic trait $\varphi\colon\Spec(R)\to M_{\alpha}^{\rm simp}$ through $L\in \Pic_{\alpha|_C}^0(C)\subset M_\alpha^{\rm simp}$ and let $\bar\varphi\colon \Spec(R)\to |h|$ be its projection. In particular,
the generic point $\Spec(K)\to |h|$ will take image in $|h|_{\rm sm}$ and, since
$\Pic_\alpha^0(\kc/|h|_{\rm sm})\subset M_\alpha$, it comes with a natural lift $\Spec(K)\to M_{\alpha}$.
By properness of the projection $M_\alpha\to|h|$, this lift extends to a trait
$\psi\colon \Spec(R)\to M_\alpha$ lifting $\bar\varphi$. Its closed point $\psi(0)\in M_{\alpha t}\subset M_\alpha$
is necessarily a smooth point of the fibre $M_{\alpha t}$.
\end{proof}

\begin{remark}\label{rem:altanal}
Here is a sketch of an alternative argument for the key step in the above proof.\smallskip

As in  \S\! \ref{eqn:MS0}, see also Remark \ref{rem:Picstable} and \S\! \ref{sec:UntwiNeron} below, we 
let $M=M(v_0)$ be moduli  space of stable sheaves with Mukai vector $v_0=(0,h,-h^2/2)$, which we view as a compactification of $\Pic^0(\kc/|h|)$. For $h$ primitive and if stability is taken with respect to an appropriate generic polarisation, $M$ is indeed a smooth projective hyperk\"ahler
variety and the projection $M\to|h|$ comes with a section given by $C\mapsto\ko_C$,
cf.\ Remark \ref{rem:Picstable}.  For $h$ primitive, the analytic Tate--{\v{S}}afarevi{\v{c}} group $\Sha(S,h)^{\rm an}$ introduced in \cite[\S\! 4.5]{HM1},
but see also \cite{Abash,AbRo,Mark,SV}, is a quotient of the affine line $\CC$ and over the affine line there exists a smooth proper family $\km\to |h|\times\CC\to\CC$ of all twists (most of them non-algebraic) of $M\to |h|$, 
 the degenerate twistor line.\smallskip

For $z\in \CC$ mapping to $\alpha$ under the projection $\CC\twoheadrightarrow\Sha(S,h)$, the fibre $\km_z\to|h|$ is birational to $M_\alpha\to|h|$ and, therefore, provides another hyperk\"ahler model
of $\Pic^0_\alpha(\kc_\eta)$. Thus,  $M\to |h|$ and $\km_z\to|h|$ are topologically isomorphic. 
In particular, in order to show that every fibre of $\km_z\to|h|$ admits a smooth point,
it suffices to show the same statement for $M\to|h|$, which is a consequence of the existence of the section $C\mapsto \ko_C$.

\end{remark}

\begin{remark}
The previous remark suggests a link between the Tate--{\v{S}}afarevi{\v{c}} group $\Sha(S,h)$ and
the Brauer group of $M=M(v_0)=\overline\Pic^0(\kc/|h|)$. As explained above, classes in $\Sha(S,h)$
induce twists of the fibration $M\to|h|$ and, motivated by the case of elliptic K3 surfaces, one
may expect the same for classes in $\Br(M)$.\smallskip

Indeed, the Leray spectral sequence
induces a map $\Br(M)\to H^1(\eta, \Pic(M_\eta))$ with $M_\eta\cong\Pic^0(\kc_\eta)$.
Since $\Pic^0(\kc_\eta)$ is principally polarised and hence $\Pic^0(M_\eta)\cong \Pic^0(\kc_\eta)$,
the  inclusion $\Pic^0(M_\eta)\,\hookrightarrow \Pic(M_\eta)$ induces a natural map
$H^1(\eta,\Pic^0(\kc_\eta))\to H^1(\eta,\Pic(M_\eta))$. However, this map is typically not an isomorphism and $\Br(M)$ only maps naturally to the latter.\smallskip

The precise link between $\Sha(S,h)$ and $\Br(M)$ has been established by \cite[Prop.\ 5.6]{MM}.
There are two cases to be distinguished depending on the divisibility of $h$, i.e.\ the integer $m$ determined by $m\ZZ=(\NS(S).h)$. For $m$ odd one has indeed an isomorphism
$\Sha(S,h)\cong \Br(M)$ while for $m$ even either again $\Sha(S,h)\cong\Br(M)$ or there is a short exact sequence
$$\xymatrix{0\ar[r]&\ZZ/2\ZZ\ar[r]&\Sha(S,h)\ar[r]&\Br(M)\ar[r]&0.}$$
\end{remark}

\subsection{N\'eron models} Our definition of a good hyperk\"ahler model is motivated by recent results of Kim \cite{Kim}. For the reader's convenience, we briefly outline those parts of his work that are relevant for us.\smallskip

 Assume $f\colon X\to B=\PP^n$ is a Lagrangian fibration of a smooth projective hyperk\"ahler variety. Assume $B_0\subset B$ is the open complement of the discriminant locus, so that the restriction of $f$ defines a smooth projective morphism $f_0\colon X_0\to B_0$. A \emph{N\'eron model}
of $f_0$ is a smooth morphism $f^{\rm n}\colon X^{\rm n}\to B$ from a possibly non-separated algebraic space $X^{\rm n}$ with the following universal property:

$\bullet$ 
Assume $g\colon Z\to B$ is a smooth morphism from a possibly non-separated algebraic space $Z$.
Then any morphism $\varphi_0\colon Z_0\coloneqq g^{-1}(B_0) \to X_0$ over $B_0$ extends uniquely to a morphism $\varphi\colon Z\to B$ over $B$.\smallskip

Following Holmes--Molcho--Orecchia--Poiret \cite{HMOP}, Kim deviates from the standard notion of a N\'eron model which typically includes `separated',
see \cite[Rem.\ 3.2]{Kim} for a thorough discussion.
\medskip

For a good hyperk\"ahler model  $X$ of its generic fibre $A\coloneqq X_\eta$ (which is a torsor over
the abelian variety $A_0\coloneqq\Aut(X_\eta)$), 
a N\'eron model  $f^{\rm n}\colon X^{\rm n}\to B$ of $f_0$ exists. This is an application
of \cite[Thm.\ 1.2]{Kim} and uses the hypothesis (ii) in Definition \ref{def:goodHK}. Furthermore, 
\cite[Thm.\ 1.3]{Kim} implies that $f^{\rm n}$ is surjective. Indeed, Kim not only proves the existence
of a N\'eron model for $f_0$ but also of a N\'eron model $P\to B$ for its relative automorphism  scheme $P_0\coloneqq\Aut^0(X_0/B_0)\to B_0$. He shows that $P\to B$ is a smooth and commutative group algebraic space but not necessarily separated. Moreover,
$X^{\rm n}\to B$ is a torsor over $P\to B$ and, in particular, $f^{\rm n}$ is surjective. 
Note that by \cite[Thm. 1.6]{Kim}, the N\'eron model $X^{\rm n}$ contains  
$X'=X\,\setminus\,{\rm Sing}(f)$ as an open subset. However, as the same applies to other good hyperk\"ahler models
of $X_\eta$, the N\'eron model $X^{\rm n}$ is in general really not separated.\smallskip

In the following we will think of the N\'eron model $f^{\rm n}\colon X^{\rm n}\to B$ as a class
$$[X^{\rm n}]\in H^1(B,P),$$
where we identify the commutative group algebraic space $P\to B$ with its sheaf of sections.

\subsection{N\'eron model of the untwisted $\pmb{\Pic^0}$}\label{sec:UntwiNeron}
We consider again a polarised K3 surface
$(S,h)$, as always assuming $h$ primitive. Associated with $(S,h)$ there are two group spaces\footnote{The notation for the second is taken from \cite{Holmes}.}
$$\Pic^0(\kc/|h|)\to |h|~\text{ and }~\Pic^{[0]}(\kc/|h|)\to |h|.$$ The former  parametrises all
invertible sheaves on curves $C\in |h|$ that restrict to invertible sheaves of degree zero on
all integral components  $C_i\subset C$ while the latter only requires the invertible sheaves to be of total
degree zero. Here, the total degree on a possibly non-reduced curve $C=\sum m_iC_i$
is defined to be $\deg(L)\coloneqq \sum m_i\deg(L|_{C_i})$. Clearly, one is an open subset of the other:
\begin{equation}\label{eqn:Pic00}
\Pic^0(\kc/|h|)\subset \Pic^{[0]}(\kc/|h|).
\end{equation} Furthermore, if all curves in $|h|$ are integral, which is a Zariski open condition for $(S,h)$ and always satisfied when $\Pic(S)=\ZZ h$, then (\ref{eqn:Pic00}) is an equality.\smallskip

Note that, while $\Pic^0(\kc/|h|)\twoheadrightarrow |h|$ is a separated commutative
group scheme (e.g.\ because it is an open subset of the 
projective moduli space $M$, see below), the bigger $\Pic^{[0]}(\kc/|h|)$ is often not separated. Another point of view
is to consider any invertible sheaf on a curve $C\in |h|$ as a simple sheaf on the ambient K3 surface $S$. Then
$\Pic^{[0]}(\kc/|h|)$ is identified with a certain open subset of the non-separated algebraic space of all simple sheaves on $S$. 
Also note that the projection $\Pic^0(\kc/|h|)\to |h|$ is surjective and all fibres 
are connected, for the latter see \cite[\S\! 9.2, Cor.\ 13]{BLR}. The fibres of $\Pic^{[0]}(\kc/|h|)\to |h|$ are not necessarily connected anymore,
but the projection is still smooth, cf.\ \cite[\S\! 8.4, Thm.\ 1]{BLR}.
\smallskip

By abuse of notation, we write $\kc\to|h|_{\rm sm}$  for the universal family of smooth curves in the linear system $|h|$.
Then  $\Pic^0(\kc/|h|_{\rm sm})\to|h|_{\rm sm}$ is an abelian scheme with scheme-theoretic generic
fibre $\Pic^0(\kc_\eta)$, the Jacobian of the generic fibre $\kc_\eta$. It is partially compactified by
$$\Pic^0(\kc_\eta)\subset \Pic^0(\kc/|h|_{\rm sm})\subset\Pic^0(\kc/|h|)\subset\Pic^{[0]}(\kc/|h|).$$
A proper hyperk\"ahler model is provided by the moduli space of stable sheaves $$
\Pic^0(\kc_\eta)\subset\Pic^0(\kc/|h|_{\rm sm})\subset\Pic^0(\kc/|h|)\subset M,$$ see (\ref{eqn:MS0}) and Remark \ref{rem:Picstable}. However, since $\Pic^{[0]}(\kc/|h|)$ is frequently not separated, it cannot be contained in $M$.\smallskip

 By virtue of Lemma \ref{lem:GoodHK}, the results of \cite{Kim} apply
to the Lagrangian fibration $\pi\colon M\to |h|$, so that we have a N\'eron model $$M^{\rm n}\to |h|,$$
as a smooth, commutative group algebraic space, at our
disposal. In particular, $M\,\setminus\,{\rm Sing}(\pi)$ is an open subset of $M^{\rm n}$, but be warned
that $M^{\rm n}\to |h|$ is not proper and  does contain all of the moduli space $M$. Also, $M^{\rm n}$
is not separated as soon as $M\to |h|$ admits other relative minimal models.
\smallskip

The following result is intuitively clear and certainly well known to the experts at least under additional assumptions, cf.\ \cite[\S\! 9.5, Thm.\ 1]{BLR} or \cite[Prop.\ 6.10]{Kim}.

\begin{prop}\label{prop:NeronPic0}
Assume that all curves $C\in|h|$ are integral.
Then $\Pic^0(\kc/|h|)$ is a N\'eron model
of $M$ or, equivalently, of $\Pic^0(\kc/|h|_{\rm sm})$:
$$M^{\rm n}\cong \Pic^0(\kc/|h|).$$
\end{prop}

\begin{proof} Assume $Z\to |h|$ is a smooth morphism and $\varphi_0\colon
Z_0\to \Pic^0(\kc/|h|_{\rm sm})$ is a morphism relative over $|h|$. Our task is to extend
$\varphi_0$ to a morphism $\varphi\colon Z\to \Pic^0(\kc/|h|)$ relative over $|h|$.
We first explain the argument under the assumption that $M$ or, equivalently,
$\Pic^0(\kc/|h|_{\rm sm})$ is a fine moduli space (which it is not in general!).
Under this unrealistic assumption,  $\varphi_0$ corresponds to
an invertible sheaf $\kl_0$ on $Z_0\times_{|h|_{\rm sm}} \kc$. 
\smallskip

The smoothness of
$Z\to |h|$ and of the total space $\kc$ imply the smoothness of $Z\times_{|h|} \kc$, cf.\ \cite[\S\! 2.3, Prop.\ 9]{BLR}.
Thus, the invertible sheaf $\kl_0$ on the open subset $Z_0\times_{|h|_{\rm sm}}\kc$ can
be extended (a priori non-uniquely) to an invertible sheaf $\kl$ on $Z\times_{|h|} \kc$. The latter
 defines a morphism $\varphi\colon Z\to\Pic(\kc/|h|)$, the image of which is clearly
contained in $\Pic^{[0]}(\kc/|h|)$.\smallskip

If we now assume that all $C\in |h|$ are integral, then $\Pic^0(\kc/|h|)=\Pic^{[0]}(\kc/|h|)$ and this space is separated. This proves that $\varphi$ takes image in $\Pic^0(\kc/|h|)$  and that $\varphi$ as an extension of $\varphi_0$ is unique.\smallskip

Let us now explain how the above has to be modified when no universal  sheaf exists. In this case, $\varphi_0$ is viewed as a section $\varphi_0\in H^0(Z_0,R^1f_{0\ast}\GG_m)$
 and we aim at extending it to a section $\varphi\in H^0(Z,R^1f_{\ast}\GG_m)$. Here, $f_0\colon Z_0\times_{|h|_{\rm sm}}\kc\to Z_0$ and $f\colon Z\times_{|h|}\kc\to Z$ denote the first projections.
 The two spaces are compared via the following commutative diagram induced by the Leray spectral
 sequence
 $$\xymatrix{H^1(Z\times_{|h|}\kc,\GG_m)\ar@{->>}[d]\ar[r]&H^0(Z,R^1f_{\ast}\GG_m)\ar[d]\ar[r]&\Br(Z)\ar@{^(->}[d]\ar[r]&\Br(Z\times_{|h|}\kc)\ar@{^(->}[d]\\
 H^1(Z_0\times_{|h|_{\rm sm}}\kc,\GG_m)\ar[r]&H^0(Z_0,R^1f_{0\ast}\GG_m)\ar[r]&\Br(Z_0)\ar[r]&\Br(Z_0\times_{|h|_{\rm sm}}\kc).\\
 }$$
 The  surjectivity resp.\ injectivity of the indicated vertical arrows follow from the smoothness of $Z$ and $Z\times_{|h|}\kc$, see above.\smallskip
 
We denote by $\beta_0\in \Br(Z_0)$ the image of $\varphi_0$ under $H^0(Z_0,R^1f_{0\ast}\GG_m)\to\Br(Z_0)$. It is the obstruction to the existence
of an invertible sheaf on $Z_0\times_{|h|_{\rm sm}}\kc$ inducing $\varphi_0$ and is pulled back under $\varphi_0$ from
the obstruction class $\gamma\in \Br(M)\subset \Br(\Pic^0(\kc/|h|))\subset\Br(\Pic^0(\kc/|h|_{\rm sm}))$ to the existence of a universal 
invertible sheaf on $\Pic^0(\kc/|h|_{\rm sm})\times_{|h|_{\rm sm}}\kc$. If $\beta_0$ can be extended to a class $\beta\in \Br(Z)$, then a simple diagram chase proves the existence of $\varphi$.
\smallskip

For the existence of the extension $\beta$ we use the closure $\bar Z_0$ of the graph of $\varphi_0$ in $Z\times_{|h|}M$:
$$\xymatrix{Z_0\ar@{^(->}[d]\ar@{^(->}[r]&Z_0\times_{|h|_{\rm sm}}\Pic^0(\kc/|h|_{\rm sm})\ar@{^(->}[d]\\
\bar Z_0\ar@{^(->}[r]&Z\times_{|h|}M.}$$
Since $M$ is projective, the first projection $\pi\colon\bar Z_0\twoheadrightarrow Z$ is proper and birational. In particular, there is a closed subset $T\subset Z$ of codimension at least two with
$\bar Z_0\,\setminus\,\pi^{-1}(T)\cong Z\,\setminus\, T$. The second projection $\bar Z_0\to M$ allows us to pull back $\gamma\in \Br(M)$ to a class in $\Br(\bar Z_0)$. Restriction to the complement of
$\pi^{-1}(T)$ gives an extension of $\beta_0$ to a class in $\Br(Z\,\setminus\,T)$. Conclude
by using the smoothness of $Z$ which implies $\Br(Z)\cong\Br(Z\,\setminus\,T)$, see \cite[Thm.\ 6.1]{BrauerIII} or \cite[\S\! 3.7]{CTS}.
\end{proof}

\begin{remark}\label{rem:Neron[0]}
In the proof of Proposition \ref{prop:NeronPic0} we actually showed that $\Pic^{[0]}(\kc/|h|)$ always
has at least the existence property of a N\'eron model, i.e.\ any morphism $\varphi_0\colon Z_0\to \Pic^0(\kc/|h|_{\rm sm})$
can be extended to a morphism $\varphi\colon Z\to \Pic^{[0]}(\kc/|h|)$, everything relative over $|h|$. However,  $\varphi$ will not be unique in general.
\end{remark}

Also note that $\Pic^{[0]}(\kc/|h|)\to|h|$ is smooth, cf.\ \cite[\S\! 8.4, Thm.\ 1]{BLR}. Hence, by the defining property of the N\'eron model, there exists a unique
morphism 
\begin{equation}\label{eqn:Pic0Neron}
\Pic^{[0]}(\kc/|h|)\to M^{\rm n}
\end{equation}
which extends the open inclusion of $\Pic^0(\kc/|h|_{\rm sm})\subset M$. Altogether, this shows that one should think of $M^{\rm n}$ as a quotient of
$\Pic^{[0]}(\kc/|h|)$ that kills part of the non-separatedness of $\Pic^{[0]}(\kc/|h|)$
but possibly not all of it, cf.\ \cite[\S\! 6.3.3]{Kim}.

\begin{cor}\label{cor:surjH1}
For an arbitrary polarised K3 surface $(S,h)$ with $h$ primitive, the map
(\ref{eqn:Pic0Neron}) of group algebraic spaces over $|h|$ induces a surjective map
\begin{equation}\label{qn:Pic0Neron2}
H^1(|h|,\Pic^{[0]}(\kc/|h|))\twoheadrightarrow H^1(|h|,M^{\rm n}).
\end{equation}
\end{cor}

\begin{proof}
We apply Remark \ref{rem:Neron[0]} to the inclusion $\varphi_0\colon Z_0=\Pic^0(\kc/|h|_{\rm sm})\,\hookrightarrow \Pic^{[0]}(\kc/|h|)$ and  its compactification $Z=M^{\rm n}$ to obtain a (non-unique) morphism
$\varphi\colon M^{\rm n}\to \Pic^{[0]}(\kc/|h|)$ which in turn induces a map $H^1(\varphi)\colon H^1(|h|,M^{\rm n})\to H^1(|h|,\Pic^{[0]}(\kc/|h|))$ between the cohomology groups of their sheaves of sections.
\smallskip

Using the uniqueness part for the N\'eron model $M^{\rm n}$,
we can conclude that the composition of $\varphi$ with (\ref{eqn:Pic0Neron}) is the identity. Hence, $H^1(\varphi)$ is a section of (\ref{qn:Pic0Neron2}), which immediately proves the claimed surjectivity.
\end{proof}

\subsection{Surjectivity} 
As before, we consider the universal curve $p\colon\kc\to|h|$ of a polarised K3 surface $(S,h)$.
The role of the abelian variety $A_0$ (over the function field $k(\eta)$ of $|h|$) in the general setting 
of Definition \ref{def:goodHK} is here played by the scheme-theoretic generic fibre $\Pic^0(\kc/|h|)_\eta\cong \Pic^0(\kc_\eta)$.

\begin{prop}\label{prop:Surj}
Assume $(S,h)$ is a polarised K3 surface with $h$ primitive. Then the image
of the injection (\ref{eqn:inj}) 
\begin{equation}\label{eqn:im}
\Sha(S,h)\,\hookrightarrow H^1(\eta,\Pic^0(\kc_\eta)),~ \alpha\mapsto\Pic^0_\alpha(\kc_\eta)
\end{equation}
is the set of all $\Pic^0(\kc_\eta)$-torsors that admit a good hyperk\"ahler model.
\end{prop}

\begin{proof} That all torsors in the image admit a good hyperk\"ahler model is the content
of Lemma \ref{lem:GoodHK}.\smallskip

For the converse, let $A$ be a  $\Pic^0(\kc_\eta)$-torsor with a good hyperk\"ahler model
$A\subset X\to|h|$. We have to show that then $[A]\in H^1(\eta,\Pic^0(\kc_\eta))$ is contained in the image of (\ref{eqn:im}).
\smallskip

By assumption, $X_\eta=A$ is a $\Pic^0(\kc_\eta)$-torsor and, therefore,
$\Aut^0(X_\eta)$ is isomorphic to $\Pic^0(\kc_\eta)$.
(To make this more geometric, one might want to extend it to an isomorphism  
\begin{equation}\label{eqn:N1}
\Aut^0(X_U/U)\cong\Pic^0(\kc_U/U)
\end{equation} of abelian schemes over some small non-empty open
subset $U\subset|h|$.) Due to the uniqueness of the N\'eron model, (\ref{eqn:N1}) induces
an isomorphism $$P\cong M^{\rm n}$$
of the N\'eron models $P$ of $\Aut(X_0/|h|_0)$ and $M^{\rm n}$ of the moduli space
$M\to|h|$ (or, equivalently, of $\Pic^0(\kc/|h|_{\rm sm})\to|h|_{\rm sm}$). \smallskip

For better readability, we introduce the following shorthands: By $\PPic$ we denote the \'etale sheaf
of sections of $\Pic(\kc/|h|)\to |h|$, so $\PPic\cong R^1p_\ast\GG_m$, and we let $\PPic^{[0]}$
denote the \'etale sheaf of sections of $\Pic^{[0]}(\kc/|h|)\to|h|$.
Then, taking cohomology of the natural short exact sequence
\begin{equation}\label{eqn:DefA}
 \xymatrix{0\ar[r]&\PPic^{[0]}\ar[r]&\PPic \ar[r]^-{\deg} &\ZZ\ar[r]&0,}
 \end{equation}
 results in the short exact sequence
 $$\xymatrix{0\ar[r]&\ZZ/m\ZZ\ar[r]& H^1(|h|,\PPic^{[0]})\ar[r]& H^1(|h|, \PPic)\ar[r]&0.}$$
Here, we use  $H^1(|h|,\ZZ)=0$ and $\im(H^0(|h|,R^1p_\ast\GG_m)\to H^0(|h|,\ZZ))=m\ZZ$, where
 the latter follows again from the Leray spectral sequence and $\Pic(\kc)\cong\Pic(S)\oplus\ZZ p^\ast\ko(1)$. 
\smallskip

 \Old{
 To illustrate the argument, we shall first discuss the case that all curves
 in $|h|$ are integral.\TBC{I AM TRYING TO AVOID THE ANALYTIC TS GROUP as long as possible. It would feel a little unnatural for my taste.} This holds for a non-empty Zariski open subset of 
 polarised K3 surfaces $(S,h)$, a subset which contains all $(S,h)$ with $\Pic(S)\cong\ZZ h$.
 Let  $p\colon \kc\to|h|$ be the projection and  denote by $\ka$ the sheaf of sections
 of $\Pic^0(\kc/|h|)\to |h|$. By definition and under our assumption that all curves in $|h|$ be integral, $\ka$ sits in a short exact
 sequence
 \begin{equation}\label{eqn:DefA}
 \xymatrix{0\ar[r]&\ka\ar[r]&R^1p_*\GG_m \ar[r]^-{\deg} &\ZZ\ar[r]&0}
 \end{equation}
 and taking cohomology gives the short exact sequence
 $$\xymatrix{0\ar[r]&\ZZ/m\ZZ\ar[r]& H^1(|h|,\ka)\ar[r]& H^1(|h|, R^1p_\ast\GG_m)\ar[r]&0.}$$
Here, we use  $H^1(|h|,\ZZ)=0$ and $\im(H^0(|h|,R^1p_\ast\GG_m)\to H^0(|h|,\ZZ))=m\ZZ$, where
 the latter follows again from the Leray spectral sequence and $\Pic(\kc)\cong\Pic(S)\oplus\ZZ p^\ast\ko(1)$. Note that by definition,
 we have  $H^1(|h|,\ka)\cong H^1(|h|,\Pic^0(\kc/|h|))$.}

 On the other hand, the Leray spectral sequence defines
 an isomorphism $$\Br(S)\cong\Br(\kc)\cong H^1(|h|,R^1p_\ast\GG_m)=H^1(|h|,\PPic).$$
 The injectivity follows from $H^0(|h|,p_\ast\GG_m)\cong\Br(|h|)=0$ and for the surjectivity one
 uses  $H^3(|h|,p_\ast\GG_m)=H^3(|h|,\GG_m)=0$, cf.\ \cite[Prop.\ 1.4]{BrauerII}. 
Imitating the discussion in \cite[\S\! 5.2]{HM1} in the case of elliptic K3 surfaces, see also \cite{DMS},
 the picture can be completed to
 a commutative diagram
 \begin{equation}\label{eqn:CD}
 \xymatrix{0\ar[r]&\ZZ/m\ZZ\ar[r]\ar[d]^-=& \Sha(S,h)\ar[d]^-\cong\ar[r]& \Br(S)\ar[r]\ar[d]^-\cong&0\\
  0\ar[r]&\ZZ/m\ZZ\ar[r]& H^1(|h|,\PPic^{[0]})\ar[r]& H^1(|h|, \PPic)\ar[r]&0.}
  \end{equation}

By definition, $[A]$ is the image of $[X^{\rm n}]$ under $H^1(|h|,M^{\rm n})\cong H^1(|h|,P)\to H^1(\eta,\Pic^0(\kc_\eta))$
and, according to Corollary \ref{cor:surjH1}, the map 
\begin{equation}\label{eqn:onto}
\Sha(S,h)\cong H^1(|h|,\PPic^{[0]})\twoheadrightarrow H^1(|h|,M^{\rm n})\cong H^1(|h|,P)
\end{equation} is surjective. Hence, $[A]$ is contained in the image
of $\Sha(S,h)\to H^1(\eta,\Pic^0(\kc_\eta))$, which concludes the proof.
\end{proof}

As an aside, we observe that the surjection (\ref{eqn:onto}) is in fact always an isomorphism
$$\Sha(S,h)\cong  H^1(|h|,P).$$ This follows from
the injectivity of the two restriction maps $\Sha(S,h)\,\hookrightarrow H^1(\eta,\Pic^0(\kc_\eta))$, see
Corollary \ref{cor:Inj}, and $H^1(|h|,P)\,\hookrightarrow H^1(\eta, P_\eta)$, see \cite[Prop.\ 6.32]{Kim}.
Furthermore, if all $C\in |h|$ are integral, then by Proposition \ref{prop:NeronPic0} we even have isomorphisms $$P\cong M^{\rm n}\cong\Pic^{[0]}(\kc/|h|)\cong\Pic^0(\kc/|h|).$$
\section{Conclusion}
Let us briefly summarise the main steps to conclude the proof of the main result and
add a few comments on possible generalisations.
\subsection{Proof of Theorem \ref{thm:main}}\label{sec:Proof}
The map $$\Sha(S,h)\to\{\,A=\Pic^0(\kc_\eta)\text{-torsor}\,\}, ~\alpha\mapsto \Pic_\alpha^0(\kc_\eta)$$
is injective by virtue of Corollary \ref{cor:Inj} and according to Proposition \ref{prop:Surj} 
the image consists of all those torsors that admit a good hyperk\"ahler model.
This concludes the proof.\qed

\subsection{Non-primitive polarisation}
Throughout we have insisted that the polarisation $h$ be primitive.
This is the case for the motivating classical example of an elliptic K3 surface $S\to\PP^1$, in which,
however, $h$ ceases being ample. Naturally, one wonders what happens in the non-primitive case.
In our discussion, the primitivity of $h$ was used in both parts of the proof, for the injectivity in \S\! \ref{sec:injec} and for the surjectivity in \S\! \ref{sec:Image}. \smallskip

For example, in \S\! \ref{sec:injec} we needed the smooth hyperk\"ahler compactifications $M$ and $M_\alpha$ of $\Pic^0(\kc_\eta)$ resp.\ $\Pic_\alpha^0(\kc_\eta)$ and their weight two Hodge structures $H^2(M,\ZZ)$ and $H^2(M_\alpha,\ZZ)$.
If $h$ is not primitive, then smooth hyperk\"ahler compactifications usually do not exist.
Nevertheless, it is feasible that this part of the proof could be replaced by another argument not
relying on $h$ being primitive.  Also in the second part of the proof of Proposition \ref{prop:Inj} 
the primitivity of $h$ was used, namely for the existence of a twisted universal sheaf, which in turn
relies on  the moduli space parametrising only stable sheaves. Nevertheless, we believe that the injectivity should also hold without any further assumption on the ample class $h$ but both parts of the proof would need to be adapted. 
\smallskip

For the surjectivity, however, one would first need a potential characterisation of the image which
eludes us for the time being. Examples of moduli spaces with $h$ not primitive have been studied in the literature, see e.g.\  \cite[Thm.\ 1.1]{IGB} and also the discussion in \cite[\S\! 4]{dCRS}  and \cite[\S\! 6.4.5]{Kim}. In the two cases studied in \cite{dCRS} the moduli spaces are birationally smooth projective hyperk\"ahler of OG10 resp.\ ${\rm K3}^{[5]}$-type, but the structure of their fibres are more complicated. Especially, everywhere non-reduced fibres do occur, cf.\ \cite[Thm.\ 1.1]{IGB}.\smallskip
 \smallskip
 
 The assumption on $h$ was also used in Lemma \ref{lem:GoodHK} which shows that all $\Pic^0_\alpha(\kc_\eta)$  admit a good hyperk\"ahler model. This was pointed out explicitly
 in the proof, but also the approach via the analytic Tate--{\v{S}}afarevi{\v{c}} group in Remark \ref{rem:altanal} makes use of this assumption. First of all, for $h=kh_0$, $k>1$,
 $\Pic^0(\kc/|h|_{\rm sm})$
 typically does not admit a smooth hyperk\"ahler compactification. Second, if one starts with $\Pic^d(\kc/|h|_{\rm sm})$ that admits a smooth hyperk\"ahler compactification, we have no control over
 its fibres anymore. Also note that the analytic Tate--{\v{S}}afarevi{\v{c}} group in this case has more than one connected component.
 




\end{document}